\documentclass[11pt]{amsart}
\usepackage{amsfonts,eucal}
\usepackage[all]{xy}
\begin{document}

\newcommand{\A}{{\mathcal A}}
\newcommand{\Aut}{{\rm Aut}}
\newcommand{\alg}{{\rm Alg}}
\newcommand{\B}{{\mathcal B}}
\newcommand{\C}{{\mathbb C}}
\newcommand{\Cat}{{\rm Cat}}
\newcommand{\Comod}{{\rm Comod}}
\newcommand{\ex}{{\rm ex}}                   
\newcommand{\F}{{\mathbb F}}                            
\newcommand{\Fun}{{\rm Fun}}
\newcommand{\G}{{\mathbb G}}
\newcommand{\bH}{{\mathbb H}}      
\newcommand{\Hom}{{\rm Hom}}
\newcommand{\K}{{\rm K}}
\newcommand{\tK}{{{\rm K}_\dagger}}
\newcommand{\lin}{{\rm lin}}
\newcommand{\Mod}{{\rm Mod}}
\newcommand{\Mot}{{\sf Mot}}
\newcommand{\oh}{{\bf o}}
\newcommand{\op}{{\rm op}}
\newcommand{\tr}{{\rm tr \;}}
\newcommand{\perf}{{\rm perf}}
\newcommand{\pt}{{\rm pt}}         
\newcommand{\Q}{{\mathbb Q}}
\newcommand{\T}{{\mathbb T}}
\newcommand{\TC}{{\rm TC}}
\newcommand{\tC}{{{\rm TC}_\dagger}}
\newcommand{\THH}{{\rm THH}}
\newcommand{\TR}{{\rm TR}}
\newcommand{\Vect}{{\rm Vect}}
\newcommand{\Z}{{\mathbb Z}}

\title{The cosmic Galois group as Koszul dual to Waldhausen's $A(*)$}
\author{Jack Morava}
\address{Department of Mathematics, Johns Hopkins University,
Baltimore, Maryland 21218}
\email{jack@math.jhu.edu}
\thanks{This work was supported by the NSF}
\date {23 August 2011}
\begin{abstract}{\noindent The world is so full of a number of things\\
I'm sure we shall all be as happy as kings.\\
\medskip

Robert Louis Stevenson, {\it A Child's Garden of Verses}}\end{abstract}

\maketitle

\section {Basic questions} \bigskip

\noindent
{\bf 1.1  Existence:} Why is there something, rather than nothing? \bigskip

\noindent
This does not seem very accessible by current methods. A more realistic goal 
may be \bigskip

\noindent
{\bf Classification:} Given that there's something, what could it be?
\bigskip

\noindent
This suggests a \bigskip

\noindent
{\bf Program:} If things fall into {\bf categories} ($\A,\B$,\dots), hopefully {\bf small} and
{\bf stable} enough to be manageable, techniques from K-theory may be useful. \bigskip

\noindent
{\bf 1.2} In fact Blumberg, Gepner, and Tabuada ([4], see also [10]) have constructed a Cartesian 
closed category $\Cat_\infty^\perf$ of small stable $\infty$-categories, eg $\A,\B, \Fun^\ex(\A,\B), 
\dots$ and there is then a (similarly Cartesian closed) big {\bf spectral} category of {\bf pre}-motives: 
with objects as above, and morphism objects
\[
\Hom_\Mot(\A,\B) := \K (\Fun^\ex(\A,\B)) \in \K(\$) - \Mod
\]                                                               
enriched over Waldhausen's $A$-theory spectrum. [The superscript `ex' signifies functors which
preserve finite limits and colimits, and the objects of the category are taken to be idempotent 
complete (ie, the category is suitably localized with respect to Morita equivalence).] 

\newpage

\noindent
Such a category has a functorial completion to a {\bf pre} - triangulated category $\Mot$
([6 \S 4.5]: ie, whose homotopy category is triangulated); this involves enlarging the set of objects
by adjoining suitable cofibers, generalizing the classical Karoubification in Grothendieck's original 
construction of a category of pure motives. \bigskip

\noindent
{\bf 1.3} Such `big' categories allow comparisons between objects from quite different areas of
mathematics (eg homotopy theory and algebraic geometry), and they raise a host of questions.\bigskip

\noindent
This posting summarizes a talk at the Hamburg 2011 conference on structured ring spectra
\[
{\tt http://www.math.uni-hamburg.de/home/richter/hh2011.html} \;.
\]
It is concerned with the motivic (Tannakian? Galois? descent?) groups of such categories 
as a tool for sorting out their relations. It is a report on work in progress with {\bf Andrew Blumberg} and
{\bf Kathryn Hess}, without whose support it would not be even a fantasy. I also want to thank
Michael Ching, Ralph Cohen, Bjorn Dundas, and Bill Dwyer for their help, and in particular for
enduring more than their share of foolish questions. Finally, much of this work is based on ideas
of Andrew Baker and Birgit Richter, and I owe them thanks for interesting conversations over many
years, and in particular for putting together this remarkable meeting. \bigskip

\section{Some examples} \bigskip

\noindent
(of things that live in this big world of motives): \bigskip

\noindent
{\bf 2.1} If $X$ is an algebraic variety over a field $k$, and $\A_X = D^\perf(\oh_X)$ is the
derived category of quasicoherent sheaves of $\oh_X$ - modules, then the class of $\A_X$ is
a version of the classical motive of $X$. The subcategory generated by such things has 
Hom-objects naturally enriched over $\K(k-\Mod)$; a cycle map associates to a
subvariety $Z$ of $X \times Y$, a resolution of its defining sheaf ${\mathcal I}_Z$ of functions,
and thus a bimodule morphism from $X$ to $Y$ \dots \bigskip

\noindent
{\bf 2.2} This example fits in the general framework of $\mathbb{A}^1$ - homotopy theory, but 
over more general rings the subject is in flux. If $X$ is an {\bf arithmetic} variety, eg over the
spectrum of integers of a number field, Deligne and Goncharov [9] have constructed a good
category of {\bf mixed Tate} motives over Spec $\Z$, with Hom objects enriched over $\K(\Z)
\otimes \Q$. The {\bf periods} of algebraic varieties [12] define similar categories
of motives. \bigskip

\noindent
{\bf 2.3} There is a great deal of interest in {\bf noncommutative} motives over a field, perhaps
also represented by suitable derived categories of perfect objects [1] \dots 

\newpage

\noindent
but my concern in this talk is to ask how the most classical example of all,\bigskip

\noindent
{\bf 2.4 topological spaces} \bigskip

\noindent
might fit in this framework. In particular, in this new world of big motives, how does the
`underlying space' or `Betti' functor
\[
X \in {\rm Varieties \; over \; \Z \mapsto X(\C) \in Spaces}
\] 
behave? This reality check is the principal motivation for this talk. \bigskip

\section{Fiber functors and their motivic automorphism groups} \bigskip

\noindent
{\bf 3.1} There are dual approaches [3,5,11] to the study of spaces in this context, both 
involving categories of modules over ring-spectra: 
\[
X \; \mapsto \$[\Omega X_+] = FX \in A_\infty - {\rm algebras}, 
\]
and 
\[
X \mapsto [X_+,\$] = DX \; (= {\rm Spanier-Whitehead \; dual}) \in E_\infty - {\rm algebras}.
\]
The first leads to Waldhausen's $A(X) = K(\$[\Omega X_+] - \Mod)$, while the second leads to
Williams' [20] $\forall (X) = K(DX - \Mod)$; together these constructions generalize Grothendieck's 
classical covariant and contravariant versions of K-theory. \bigskip

\noindent
Both $DX$ and $FX$ are {\bf supplemented} $\$$-algebras, and in good cases (ie if $X$ is both
finite and simply-connected) then 
\[
FX \cong \Hom_\$(DX,DX), \; DX \cong \Hom_\$(FX,FY) 
\]
expresses a kind of `double centralizer' duality. \bigskip
 
\noindent                                                      
{\bf 3.2} Here I'll work with the second of these alternatives, in the category with finite
$CW$-spaces $X,Y$ as objects, and morphisms
\[
\Hom_\Mot(X,Y) \sim  \K(DX \wedge DY^\op-\Mod)
\]
defined by the K-theory spectra of {\bf right-compact} $DX - DY^\op$ - bimodules [4
\S 2.16]. This category can then be made pre-triangulated, as above. \bigskip

\noindent
There are many technical variants of this construction: for example, BGT consider both
Karoubi-Villamayor and Bass-Thomason K-theory. Later we will want to modify categories of 
this sort by completing their morphism objects in various ways, and eventually
we will be interested in constructions based on THH and its relatives (TR, TC, \dots); 
then I'll label the resulting categories by the functors defining their morphism objects. 

\newpage

\noindent
For example, the cyclotomic trace defines a monoidal spectral functor
\[
\Mot_\K \to \Mot_\TC
\]
of pre-triangulated categories (and hence a triangulated functor between their homotopy 
categories). \bigskip

\noindent
{\bf 3.3} Tannakian analogs of Galois groups are a central topic in the usual theory of motives:
complicated categories can sometimes be identified, via some kind of descent, with categories 
of representations of groups of automorphisms of interesting forgetful (monoidal, `fiber') functors
to simpler categories. Weil cohomologies (Hodge, \'etale, crystalline) are classical examples, but
the following example may be more familiar here:\bigskip

\noindent
Ordinary cohomology (with coefficients in $\F_2$ and the grading neglected), viewed as a
monoidal functor
\[
H : (\rm Spectra) \ni X \mapsto H^*(X,\F_2) \in (\F_2 - \Vect) \;,
\]
has a group-valued functor
\[
\Aut^H_\otimes : (\F_2 - \alg) \ni A \mapsto \Aut^A_\otimes(H^*(-,A))
\]
of natural automorphisms, which is (co)represented by the dual Steenrod algebra:
\[
\Aut^A_\otimes(H^*(-,A)) \; \cong \; \Hom_\alg(\A^*,A) \;.
\]
The vector-space valued functor $H^*$ thus {\bf lifts} to a functor taking values in
representations of a proalgebraic groupscheme, or (in more familiar language), in the category of
$\A^*$-comodules. \bigskip

\noindent
Here I want to look at (pre-triangulated, spectral, monoidal) categories built by reducing the
morphism objects in BGT-style categories modulo the kernel of the Dennis trace $K(\$) \to \$$
(much as we can consider the category obtained from chain complexes over $\Z$ by reducing
their internal Hom-objects modulo $p$). \bigskip

\noindent
{\bf 3.4} Hess's theory of homotopical descent [14] provides us with the needed technology: a
cofibrant replacement 
\[
\xymatrix{
K(\$) \ar[dr]^\tau \ar[rr]^{\tr} & {} & \$ \\
{} & Q(\$) \ar[ur]^\rho }
\]
(of the sphere spectrum $\$$ as $K(\$)$-algebra\begin{footnote}{Note that $\K(\Z)$ is {\bf not}
similarly supplemented over $\Z$!}\end{footnote}, with $\tau$ a cofibration, and $\rho$ a weak 
equivalence) associates a `Hessian' {\bf co-ring} spectrum
\[
Q(\$) \wedge_{\K(\$)} Q(\$) \; (= \; \THH_{\K(\$)}(\$) \; )
\]
(analogous to a Hopf-Galois object in the sense of Rognes [17]) to the Dennis trace.\bigskip

\noindent
Similarly,                              
\[
\$ \to H\F_2
\]
produces the dual Steenrod algebra
\[
Q(H\F_2) \wedge_\$ Q(H\F_2) \; \sim \; \A^* \;.
\]
The resulting theory of descent relates a $K(\$)$-module spectrum $V$ to a 
$\THH_{K(\$)}(\$) := \$_{\dagger K(\$)}$ - comodule 
\[
V_{\dagger K(\$)} \; := \;  Q(\$) \wedge_{K(\$)} V \; = \; \THH_{K(\$)}(\$,V) \;,
\] 
and 
\[
\K(DX \wedge DY^\op) \to \K(DX \wedge DY^\op)_{\dagger K(\$)}  := \K_\dagger (DX \wedge
DY^\op)
\]
defines a monoidal functor
\[
\omega_\tK : \Mot_\K \to \Mot_\tK \;,
\]
the latter category being enriched over spectra with an $\$_{\dagger K(\$)}$ - comodule action
(the analog of representations of $\Aut(\omega_\tK)$). \bigskip

\noindent
We expect a more careful version of this construction to provide {\bf effective} homotopical 
descent for a version of $\Mot_\$$ with suitably completed morphism objects [14 \S 4, \S 5.5]. \bigskip

\noindent
{\bf 3.5} The notation above is unsatisfactory: it reflects similar difficulties with notation 
for Koszul duality. In the classical case of a morphism $A \to B$ of algebras over a field $k$, 
the covariant functor
\[
V \mapsto V \otimes^L_A B := V_{\dagger B} : D(A - \Mod) \to D(A_{\dagger B}-\Comod)
\]
has a contravariant $k$-vector-space dual
\[                                                          
V \mapsto V^\dagger_B := (V_{\dagger B})^* \cong {\rm RHom_A}(V,B) 
\]
with values in some derived category of ${\rm RHom}_A(B,B) := A^\dagger_B$-modules [Cartan-Eilenberg 
VI \S 5], which is in good cases a (Koszul) duality. In the formulation above, 
\[
\$_{\dagger K(\$)} \;:= \; \$ \otimes^L_{K(\$)} \$ \; = \; \THH_{\K(\$)}(\$)
\]
is the analog of the algebra of functions on a group object, while
\[
\$^\dagger_{K\$)} \; := \; {\rm RHom}_{K(\$)}(\$,\$) 
\]
is the analog of its (convolution, $L^1$) group algebra. \bigskip

\section{Cyclotomic variants} \bigskip

\noindent
{\bf 4.1} The constructions above have a straightforward analog 
\[
\Mot_\TC \to \Mot_\tC 
\]
built from topological cyclic homology; where now
\[
\tC( - ) := \THH_{\TC(\$)}(\$,\TC( - )) \in \THH_{\TC(\$)}(\$) := \$_{\dagger\TC} - \Comod 
\]
(with profinite completions implicit but suppressed)\begin{footnote}{Another interesting variant 
can be built from THH, regarded as a $\T$-equivariant spectrum.}\end{footnote}.\bigskip

\noindent
The cyclotomic trace
\[
\K(\$) \to \TC(\$) \sim \$ \vee \Sigma \C P^\infty_{-1}
\]
(again mod completion) identifies the K-theory spectrum with $\$ \vee \Sigma \bH P^\infty_+$ at 
regular odd primes [15, 18]. The cofibration
\[
S^{-1} \to \Sigma \C P^\infty_{-1} \to \Sigma \C P^\infty_+ 
\]
suggests that the Koszul dual of $\THH_{\TC(\$)} \$$ should be close to the tensor $\$$ - algebra
$\$[\Omega \Sigma \C P^\infty_+]$ on $\C P^\infty_+$ [2]. In any case, $\$_{\dagger K(\$)}
\otimes \Q$ can be identified with the algebra of quasisymmetric functions over $\Q$, ie the
algebra of functions on a pro-unipotent group with free Lie algebra. The cyclic structure on THH
endows this Lie algebra with a $\T$-action and thus with a grading, placing one generator in each
odd degree [7]. \bigskip
 
\noindent
This is very similar to Deligne's motivic group for the category of mixed Tate motives, itself
modeled on Shafarevich's conjectured description of the absolute Galois group of $\Q$ as a
profree profinite extension of $\hat{\Z}^\times$. It leads to the appearance of odd zeta-values
in differential topology, systematically parallel to the appearance of {\bf even} zeta-values
(ie, Bernoulli numbers) in homotopy theory. \bigskip
 
\noindent
{\bf 4.2} One concern with these constructions is that neither K nor TC is {\bf
linear}, in the sense of the calculus of functors. \bigskip

\noindent
$\THH_\$(DX)$ is the realization of a cyclic object
\[
n \mapsto (DX)^{\wedge(n+1)} \; \sim \; D(X^{n+1}) 
\]
$S$-dual to the totalization of a (cocyclic) cosimplicial space modelling the free loopspace
$LX$ (cf [13]; thanks to WD for this reference!). My hope is that the homotopy fixed points
$\THH_\$(DX)^{h\T}$ can be identified as something like
\[
[E\T_+,[LX_+,\$]]^{h\T} = [LX_{h\T+},\$] = [LX_+,[E\T_+,\$]]^{h\T}
\]
and that consequently $TC(DX)$ will be accessible as a homotopy limit of things like
$[LX_+,\THH_\$(\$)]^{C_n}$. \bigskip

\noindent
This suggests that the inclusion $X \to LX$ of fixed points defines a kind of coassembly
[20] map
\[
\TC(DX) \to [X_+,\TC(\$)] 
\]
as a $\TC(\rm holim) \to {\rm holim}(\TC)$ interchange. [The classical assembly map defines a 
composition
\[
\Hom_{K(\$)}(K(\$[\Omega X_+],\$) \to \Hom_{K(\S)}(X \wedge K(\$),\$) \sim DX \; \dots]
\]
\bigskip

\noindent
{\bf 4.3} If this is so, then we can add a third step
\[
\Mot_\TC \to \Mot_\tC \to \Mot^\lin_\tC
\]
to the sequence of pre-triangulated monoidal functors above, with 
\[
\Hom^\lin_\tC(X,Y) = \THH_{\TC(\$)}(\$,[DX \wedge DY^\op,\TC(\$)]) \in 
\$_{\dagger\TC} - \Comod \;.
\]
Note that 
\[
\Hom^\lin_\tC(X,Y) \otimes \Q = HH_{\TC_\Q(\$)}(\TC_\Q(\$),H^*(Y \wedge DX))
\]
\[
= H^*(Y \wedge DX,\Q) = [Y,X]_\Q \;,
\]
so the rationalization of $\Mot^\lin_\tC$ reduces to the (rationalized) category of finite spectra, 
(conjecturally!) reconciling the motive of an algebraic variety with the stable homotopy type
of its underlying space. More generally, 
\[
[X,K(\$)]_{\dagger K(\$)} \; \sim \; [X,\$] \dots
\]
\bigskip

\bibliographystyle{amsplain}

\end{document}